\documentclass[11pt]{article}
\usepackage{dsfont}
\usepackage{graphics}
\usepackage[demo]{graphicx}
\usepackage{epsfig}
\usepackage{pstricks}
\usepackage[normal]{subfigure}
\usepackage[latin1]{inputenc}
\usepackage[english,francais]{babel}
\usepackage{relsize,exscale}
\usepackage{makeidx}
\usepackage{enumitem}
\usepackage{amsfonts,amssymb,amsmath}
\usepackage{graphicx}
\usepackage{color}
\usepackage{multirow}
\usepackage{mathrsfs}
\usepackage[normalem]{ulem}
\usepackage{cancel}
\newenvironment{prooff}{{\it Proof :}}{\hfill\rule{2mm}{2mm}\vskip3mm\par}
\newtheorem{theorem}{Theorem}[section]
\newtheorem{lemma}[theorem]{Lemma}

\newtheorem{e-definition}[theorem]{Definition\rm}
\newtheorem{remark}{\it Remark\/}

\usepackage{color}
\definecolor{dred}{rgb}{0.92,0,0}
\definecolor{dgreen}{rgb}{0,0.92,0}
\definecolor{dblue}{rgb}{0,0,0.92}
\definecolor{dyellow}{rgb}{0.95,0.95,0}

\newcommand{\R}{\mathbb{R}}
\newcommand{\N}{\mathbb{N}}
\newcommand{\A}{\mathbb{A}}
\def\D{\displaystyle}
\newcommand{\hs}{\hspace{0.1cm}}
\newcommand{\vs}{\vspace{0.4cm}}
\newcommand{\sa}{\\ [0.2cm]}
\usepackage{multirow}
\graphicspath{
{./Figures/}
{./}
}
\title{Explicit $k-$dependence for  $P_k$ finite elements in $W^{m,p}$ error estimates: application to probabilistic laws for accuracy analysis}
\author{Jo\"el Chaskalovic \thanks{D'Alembert,
Sorbonne University, Paris, France, (\emph{email}: jch1826@gmail.com)}
\qquad
Franck Assous
\thanks{
Ariel University, 40700 Ariel, Isra\"el, (\emph{email}: franckassous55@gmail.com).}
}
\date{}
\begin{document}
\maketitle
\selectlanguage{english}
\begin{abstract}
\noindent We derive an explicit $k-$dependence in $W^{m,p}$ error estimates for $P_k$ Lagrange finite elements. Two laws of probability are established to measure the relative accuracy between $P_{k_1}$ and $P_{k_2}$ finite elements ($k_1 < k_2$) in terms of $W^{m,p}$-norms. We further prove a weak asymptotic relation in $D'(\R)$ between these probabilistic laws when difference $k_2-k_1$ goes to infinity. Moreover, as expected, one finds that $P_{k_2}$ finite element is {\em surely more accurate} than $P_{k_1}$, for sufficiently small values of the mesh size $h$. Nevertheless, our results also highlight cases where  $P_{k_1}$ is {\em more likely accurate} than $P_{k_2}$, for a range of values of $h$. Hence, this approach brings a new perspective on how to compare two finite elements, which is not limited to the rate of convergence.
\\[0.2cm]
\noindent{\footnotesize {\em keywords}: Error estimates, Finite elements, C\'ea Lemma, Bramble-Hilbert lemma, Banach Sobolev spaces, Probabilistic laws.}
\end{abstract}
\section{Introduction}\label{intro}
\noindent Recently (\cite{ChasAs18} and \cite{Arxiv2}), we proposed new perspectives on relative finite elements accuracy based on a mixed geometrical-probabilistic interpretation of the error estimate derived from Bramble-Hilbert lemma. \sa
This led us to derive two laws of probability that estimate the relative accuracy, considered as a random variable, between two finite elements $P_{k_1}$ and $P_{k_2}$ ($k_1 < k_2$).\sa
%
%
By doing so, we obtained new insights which showed, among others, which of $P_{k_1}$ or $P_{k_2}$ is the most likely accurate, depending on the value of the mesh size $h$ which is no more considered as going to zero, as in the usual point of view.\sa
These results have been obtained by considering a second-order elliptic variational problem set in the Sobolev space $H^1(\Omega)$.
However, many partial differential equations are well posed in a more general class of Sobolev spaces, namely, $W^{m,p}(\Omega), (m,p)\in\N^{*2}$.\sa
Possible applications for studying case $p \neq 2$ can be the Laplace equation set in an open-bounded domain $\Omega \subset \R^n$ with a given right-hand side $f \in L^{p}(\Omega), (p\ne 2)$. Indeed, in that case, the solution to the associated variational formulation, $u$, belongs to $W^{1,p}(\Omega)$ for $p\ne2$ if the domain $\Omega$ is regular enough: this problem is indeed discussed in \cite{Brezis} (note in the Chapter on Sobolev spaces, where a reference to  \cite{AgDn59} is quoted). Other examples may be found for instance in  \cite{Gris92},  \cite{Haubxx}, \cite{KuMi12}, or in \cite{Lion69} for non-linear problems.\sa
Here, we consider  a functional framework defined by the help of $W^{m,p}$ Sobolev spaces, particularly when $p\ne 2$, and extend our previous work \cite{Arxiv2} limited to the case of the $H^1$ Hilbert space.\sa
The paper is organized as follows. We recall in Section \ref{Second_Order_Elliptic} the mathematical problem we consider as well as the basic definitions of functional tools we will need along the paper. In Section \ref{Pk_Properties}, we introduce $P_k(K)$, the space of polynomial functions defined on a given $n$-simplex $K$, of degree less than or equal to $k$. We then obtain several estimates to upper-bound the basis functions of $P_k(K)$ and their partial derivatives. We provide in Section \ref{explicit_estimate} results that make explicit the dependence of the constant involved in the \emph{a priori} $W^{m,p}-$error estimates with respect to degree $k$ of the concerned $P_k$ Lagrange finite element. Section \ref{FEM_accuracy} presents applications to the analysis of the relative finite elements accuracy in $W^{m,p}$. In particular, extending to $W^{m,p}$ spaces the two generalized probabilistic laws introduced in \cite{ChasAs18}, we prove, relying on distributions theory, and under some {\em ad hoc} assumptions that are fulfilled in many cases, that an asymptotic relation exist between these two laws. Concluding remarks follow.

\section{The abstract problem}\label{Second_Order_Elliptic}
\noindent In this section, we introduce the abstract framework we will use to derive error estimates in the $W^{m,p}$  Sobolev spaces, particularly in the non standard cases $p \neq 2$, corresponding to non-Hilbert spaces. As a consequence, we need a well-posedness result based on a stability (or inf-sup) condition extended to non-Hilbert spaces. For the error analysis, we will also need an extension of C\'ea's Lemma to Banach spaces, devoted to the approximation of  the abstract problem using a Galerkin method.\sa
\noindent In order to provide sufficient resources for a reader even not familiar with these methods to understand the approach as a whole, we recall here some fundamental results. To this end, we basically follow the presentation and the terminology proposed in the book by A. Ern and J. L. Guermond \cite{Ern_Guermond}. The book of Brenner et al. \cite{BrSc08}, that goes back to a paper by Rannacher and Scott \cite{RaSc82} can also provide helpful references. A well-informed reader may skip to subsection \ref{example}.

\subsection{Preliminary results}\label{basicBanach}
\noindent Let $W$ and $V$ be two Banach spaces equipped with their norms $\|.\|_W$ and $\|.\|_V$,  respectively. In addition, $V$ is assumed to be reflexive. Let $u \in W$ be the solution to the variational formulation
\begin{equation}\label{VP}
\left\{
\begin{array}{l}
\mbox{Find } u \in   W \mbox{ solution to:} \\ [0.1cm]
a(u,v) = l(v), \quad\forall v \in V,
\end{array}
\right.
\end{equation}
where $l$ is a continuous linear form on $V$, and $a$ is a continuous bilinear form on $W \times V$, i.e.
$$
\forall (u,v)\in W\times V,\, |a(u,v)|\leq\|a\|_{W,V}\|u\|_W\|v\|_V,
$$
with $\D\|a\|_{W,V}\equiv\inf\left\{C\in\R^{*}_{+},\forall (u,v)\in W\times V: |a(u,v)|\leq C\|u\|_W\|v\|_V\right\}$.\label{Norme_a}
Assuming that
\begin{description}
\item[\textbf{(BNB1)}] \hspace{4cm}$\D \exists \alpha > 0, \hs \inf_{w\in W}\sup_{v\in V}\frac{a(w,v)}{\|w\|_{W}\|v\|_{V}} \geq \alpha $,\vspace{0.2cm}
\item[\textbf{(BNB2)}] \hspace{3cm}$\forall v\in V, (\forall w\in W, a(w,v)=0)\Longrightarrow (v=0)$\,,
\end{description}
one can prove that problem (\ref{VP}) has one and only one solution in $W$, (see \cite{Ern_Guermond} Theorem 2.6), where \textbf{(BNB1)-(BNB2)} refers to the Banach-Necas-Babuska conditions. \sa
Now, let us introduce the approximation $u_{h}$ of $u$, solution to the approximate variational formulation
\begin{equation}\label{VP_h}
\left\{
\begin{array}{l}
\mbox{Find } u_{h} \in   W_h \mbox{ solution to:} \\
a(u_{h},v_{h}) = l(v_{h}),\quad \forall v_{h} \in V_h,
\end{array}
\right.
\end{equation}
where $W_h \subset W$ and $V_h \subset V$ are two finite-dimensional subspaces of $W$ and $V$. As noted in \cite{Ern_Guermond}, (Remark 2.23, p.92), neither condition \textbf{(BNB1)} nor condition \textbf{(BNB2)} implies its discrete counterpart. The well-posedness of (\ref{VP_h}) is thus equivalent to the two following discrete conditions:
\begin{description}
\item[(BNB1$_h$)] \hspace{4cm}$\D \exists \alpha_h > 0, \hs \inf_{w_h\in W_h}\sup_{v_h\in V_h}\frac{a(w_h,v_h)}{\|w_h\|_{W_h}\|v_h\|_{V_h}} \geq \alpha_h$,\vspace{0.2cm}
\item[(BNB2$_h$)] \hspace{3cm}$\forall v_h\in V_h, (\forall w_h\in W_h, a(w_h,v_h)=0)\Longrightarrow (v_h=0)$.
\end{description}
\vs
From now on, we assume hypotheses \textbf{(BNB1)-(BNB2)} and \textbf{(BNB1$_h$)-(BNB2$_h$)} which guarantee the well-posedness of (\ref{VP}) and (\ref{VP_h}).\sa
The last key ingredient we need for the error estimates is the following generalized C\'ea's Lemma \cite{Ern_Guermond} valid in Banach spaces:
\begin{lemma}\label{Lemme_Cea}
\textbf{(C\'ea).} Assume that $V_h\subset V$, $W_h\subset W$ and dim$(W_h)$ = dim$(V_h)$. Let $u$ solve the problem (\ref{VP}) and $u_h$ the problem (\ref{VP_h}). Then, the following error estimate holds:
\begin{equation}\label{Cea_Banach}
\D\|u-u_h\|_{W} \leq \left(1+\frac{\|a\|_{W,V}}{\alpha_h}\right)\inf_{w_h\in W_h}\|u-w_h\|_{W}.
\end{equation}
\end{lemma}
\noindent In the rest of this paper, we will consider the variational formulation (\ref{VP}) and its approximation (\ref{VP_h}) in the case where the Banach space $W$ and the reflexive Banach space $V$ are chosen as
\begin{equation}\label{V_and_W}
W\equiv W^{m,p}(\Omega) \mbox{ and } V\equiv W^{m',p'}(\Omega)\,.
\end{equation}
Above, $m$ and $m'$ are two non zero integers, $p$ and $p'$ two real positive numbers satisfying $p\ne 2$ and $p'>1$ with
\begin{equation}\label{Conjugated}
\D \frac{1}{p}+\frac{1}{p'}=1.
\end{equation}
As usual, for any integer $m$ and $1 < p < +\infty$, $W^{m,p}(\Omega)$ denotes the Sobolev space of (class of) real-valued functions which, together with all their partial distributional derivatives of order less or equal to $m$, belongs to $L^p(\Omega)$:
$$
\D W^{m,p}(\Omega) = \left\{\!\!\frac{}{}u \in L^p(\Omega)\,/\,\forall\, \alpha, |\alpha|\leq m, \partial^{\alpha}u\in L^p(\Omega)\right\},
$$
$\alpha=(\alpha_1, \alpha_2, \ldots, \alpha_n) \in \N^{n}$ being a multi-index whose length $|\alpha|$ is given by $|\alpha|=\alpha_1+\dots+\alpha_n$, and $\partial^{\alpha}u$ the partial derivative of order $|\alpha|$ defined by:
$$
\D \partial^{\alpha}u \equiv \frac{\partial^{|\alpha|}u}{\partial x_{1}^{\alpha_1}\dots\partial x_{n}^{\alpha_n}}.
$$
The norm $\|.\|_{m,p,\Omega}$ and the semi-norms $|.|_{l,p,\Omega}$ are respectively defined by:
$$
\D \forall u \in\,W^{m,p}(\Omega): \|u\|_{m,p,\Omega} = \left(\sum_{|\alpha|\leq m}\|\partial^{\alpha}u\|^{p}_{L^p}\right)^{1/p}, \hs\hs|u|_{l,p,\Omega} = \left(\sum_{|\alpha|= l}\|\partial^{\alpha}u\|^{p}_{L^p}\right)^{1/p}, 0 \leq l \leq m,
$$
where $\|.\|_{L^p}$ denotes the standard norm in $L^p(\Omega)$. \sa

\subsection{A simple example}\label{example}
\noindent We illustrate below, through an elementary example, the choice of the spaces $W$ and $V$ defined by (\ref{V_and_W}).\sa
Let $f$ be a given function that belongs to $L^{p}(]0,1[), (p\ne 2),$ and $u\in W^{2,p}(]0,1[)$ solution to:
$$
\left\{
\begin{array}{l}
-u''(x) + u(x) = f(x), x \in ]0,1[, \\[0.1cm]
u(0)=u(1)=0.
\end{array}
\right.
$$
The corresponding variational formulation is given by:
\begin{equation}\label{VP_0}
\left\{
\begin{array}{l}
\mbox{Find } u \in   W^{1,p}_{0}(]0,1[), \mbox{ solution to:} \\ [0.1cm]
\D \int_{0}^{1}\left[ u'(x)v'(x)+u(x)v(x)\right] dx = \int_{0}^{1}f(x)v(x) \,dx , \forall v \in W^{1,p'}_{0}(]0,1[),
\end{array}
\right.
\end{equation}
where $p$ and $p'$ satisfy (\ref{Conjugated}), and $W^{1,p}_{0}(]0,1[)$ denotes the space of functions $w$ of $W^{1,p}(]0,1[)$ such that $w(0)=w(1)=0$.
\begin{remark}$\frac{}{}$
\begin{itemize}
\item First of all, we notice that all the integrals in (\ref{VP_0}) are bounded due to H\''older's inequality.
\item Second, taking for example $p=3/2$ and $q=3$, the corresponding spaces $W$ and $V$ introduced above are equal to
$W=W^{1,3/2}_{0}(]0,1[)$ and $V=W^{1,3}_{0}(]0,1[)$, that are respectively a Banach space and a reflexive Banach space, as required.
\end{itemize}
\end{remark}
\noindent In the rest of the paper, we shall assume that $\Omega$ is an open subset in $\R^n$, exactly covered by a mesh ${\mathcal T}_h$ composed by $N_K$ $n$-simplexes $K_{\mu}, (1\leq \mu \leq N_K),$ which respect classical rules of regular discretization, (see for example \cite{ChaskaPDE} for the bidimensional case, or \cite{RaTho82} in $\R^n$). Moreover, we denote by $P_k(K_{\mu})$ the space of polynomial functions defined on a given $n$-simplex $K_{\mu}$ of degree less than or equal to $k$, ($k \geq$ 1). \sa
Henceforth, we assume that the approximate spaces $W_h$ and $V_h$, satisfying dim$(W_h)$ = dim$(V_h)$, are included in the space of functions defined on $\Omega$, composed by polynomials belonging to $P_k(K_{\mu}), (1 \leq \mu \leq N_K)$. As a consequence, $W_h \subset W^{m,p}(\Omega)$ and $V_h \subset W^{m',p'}(\Omega)$.\sa
In the following section, we derive appropriate estimates related to the canonical basis of $P_k(K_\mu)$. This will in turn enable us to make explicit the dependence on $k$ of the constant involved in the \emph{a priori} error estimates in $W^{m,p}(\Omega)$.
\section{Properties of Lagrange finite element $P_k$}\label{Pk_Properties}
\noindent In this section we follow the definitions and properties of the $P_k$ finite element in $\R^n$ described by P. A. Raviart and J. M. Thomas in \cite{RaTho82}. \sa
Let us consider a $n$-simplex $K \subset \R^n$ which belongs to a regular mesh ${\mathcal T}_h$. Since a complete polynomial of order $k$ which belongs to $P_k(K)$ contains
\begin{equation}\label{Dim_Pk}
\D N \equiv \left(
  \begin{array}{c}
n+k \\
n
  \end{array}
\right)
= \frac{(n+k)!}{n!\,k!}
\end{equation}
terms, each $n$-simplex element of the mesh ${\mathcal T}_h$ must be associated with $N$ independent specifiable parameters, or degrees of freedom, to assure the unisolvence of the finite element \cite{RaTho82}. \sa
Then, it is convenient to carry out all analysis of $n$-simplexes in terms of the so-called $n$-simplex barycentric coordinates $\lambda_1,\dots,\lambda_{n+1}$ which satisfy $\D \sum_{i=1}^{n+1}\lambda_i=1$.\sa
A regularly spaced set of points $M_{i_1,\dots,i_{n+1}}$ can be defined in a $n$-simplex $K$ by the barycentric coordinates values
$\D M_{i_1,\dots,i_{n+1}} = \left(\frac{i_1}{k},\dots,\frac{i_{n+1}}{k}\right), \hs 0 \leq i_1,\dots,i_{n+1} \leq k$ satisfying
\begin{equation}\label{sum_ij_egal_k}
i_1 + \dots +i_{n+1}=k.
\end{equation}
One can check that the number of points defined in this way is equal to $N$, the dimension of $P_k(K)$ given by (\ref{Dim_Pk}). \sa
Therefore, we introduce the canonical basis of functions $p_{i_1,\dots,i_{n+1}}$ of the variables $(\lambda_1, \dots, \lambda_{n+1})$ which belongs to $P_k(K)$ defined by:
\vspace{-0.3cm}
\begin{equation}\label{shape_function}
\D p_{i_1, \dots, i_{n+1}}(\lambda_1, \dots, \lambda_{n+1}) \equiv \prod_{j=1}^{n+1}P_{i_j}(\lambda_j),
\end{equation}
where the auxiliary polynomials $P_{i_j}(\lambda_j)$ are given by:
\begin{equation}\label{P_ij}
\D \D P_{i_j}(\lambda_j) \equiv \left |
\begin{array}{ll}
\hs \D\prod_{c_j=1}^{i_j}\left(\frac{k \lambda_j - c_j +1}{c_j}\right), & \mbox{ if } \hs i_j \geq 1, \vspace{0.1cm} \\
\hs 1, & \mbox{ if } \hs i_j = 0.
\end{array}
\right.
\end{equation}
$P_{i_j}$ is clearly a polynomial of order $i_j$ in $\lambda_j$, and therefore, due to condition (\ref{sum_ij_egal_k}), $p_{i_1, \dots, i_{n+1}}$ given by (\ref{shape_function}) is a polynomial of order $k$. \sa
In the sequel, we will also use a single-index numbering to substitute the multi-index one. It will be the case for the $N$ points $M_{i_1,\dots,i_{n+1}}$ simply denoted $(M_i)_{i=1,N}$, as well as for the $N$ canonical functions $p_{i_1,\dots,i_{n+1}}$ denoted $(p_i)_{i=1,N}$, and so on. \sa
Let us also remark that each polynomial $p_i$ defined by (\ref{shape_function})-(\ref{P_ij}) is characteristic to the corresponding point $M_i$. That is to say that we have the following property (see \cite{RaTho82}):
$$
\forall\, i,j=1\mbox{ to } N: p_i(M_j)=\delta_{ij}.
$$
Therefore, for a given set of $N$ values $\varphi_i \equiv \varphi_{{i_1, \dots, i_{n+1}}}$ known at the $N$ points $M_i \equiv M_{i_1, \dots, i_{n+1}}$, the polynomial $Q$ in $P_k(K)$ given by:
\begin{eqnarray*}
\D\forall M \in K: Q(M) & = & Q(\lambda_1, \dots, \lambda_{n+1}), \nonumber\\[0.1cm]
& = & \hspace{-0.7cm}\sum_{i_1 + \dots + i_{n+1}=k} \hspace{-0.5cm}\varphi_{i_1, \dots, i_{n+1}}\,p_{i_1, \dots, i_{n+1}} \! (\lambda_1, \dots, \lambda_{n+1}) \, = \, \sum_{i=1}^{N} \varphi_i p_i(\lambda_1, \dots, \lambda_{n+1}),
\end{eqnarray*}
is the unique one in $P_k(K)$ such that $Q(M_i)= \varphi_i$. \sa

\noindent The following lemma gives the first point-to-point estimates for the polynomials $p_i$ defined by (\ref{shape_function})-(\ref{P_ij}).
\begin{lemma}
Let $(p_i)_{i=1,N}$ be the canonical basis functions of the space of polynomials $P_k(K)$ which are defined by (\ref{shape_function})-(\ref{P_ij}). \sa
Then:
\begin{equation}\label{deriv_part_ordre_m_pi_vs_lambda_l}
\D|p_i(\lambda_1, \dots, \lambda_{n+1})|\leq k^{n+1}, \hs \forall\, r \in \N^*: \left|\frac{\partial^{\,r} p_i}{\partial\lambda_{q_{1}}\dots\partial\lambda_{q_{r}}}(\lambda_1, \dots, \lambda_{n+1})\right| \leq k^{r(n+2)},
\end{equation}
where $(q_{1}, q_{2}, \ldots, q_{r}) \in \N^r$.
\end{lemma}
\begin{prooff} $\frac{}{}$ \sa
This lemma generalizes lemma 5.2 of \cite{Arxiv2} in the case where $p$ is not necessarily equal to $2$. Hence, we will only provide a sketch of the proof, and refer the interested reader to this reference for details. \sa
\noindent The logical sequence of the proof can be summarized as follows:\sa
\noindent $\blacktriangleright$ First, examine the upper boundary of the basis functions $p_i, (i= 1,\dots, N)$. \sa
This requires to introduce integer $n_i$ $(0 \leq n_i\!\leq n+1)$ corresponding to the number of polynomials $P_{i_j}(\lambda_j)$ such that:
\begin{eqnarray*}
\forall j=1,\dots,n_i, \,(n_i \geq 1),\,& : & P_{i_j}(\lambda_j)=P_{1}(\lambda_j)=k\lambda_j, (i_j=1), \\[0.1cm]
\forall j=n_i+1,\dots,n+1, \,(n_i \leq n)\, & : & P_{i_j}(\lambda_j)=\frac{k\lambda_j(k\lambda_j - 1)\dots(k\lambda_j - i_j+1)}{i_j!}, (i_j>1).
\end{eqnarray*}
Using the fact that the structure of $p_i$ depends on the value of $n_i$, ($n_i=0$, $1 \leq n_i \leq n$ or $n_i=n+1$), one obtains in all cases that
$$
\D |p_i(\lambda_1,\dots,\lambda_{n+1})| \leq k^{n+1}.
$$
\noindent  $\blacktriangleright$  Consider next $r=1$,  which corresponds to upper-bound the partial derivative $\D\frac{\partial p_i}{\partial\lambda_q}$, for a given pair of non-zero integers $(i,q)$. \sa
Once again, depending on the value value of $n_i$, we obtain different estimates, the more restrictive being
$$
\D\left|\frac{\partial p_i}{\partial\lambda_q}\right| \leq k^2 \, k^{n_i} \leq k^{n+2}.
$$
\noindent $\blacktriangleright$ Finally, handle the partial derivative of $p_i$ of order $r$ with respect to $\lambda_{q_1},\dots,\lambda_{q_r}$. \sa
To this end, we basically use the upper bound and remark that any first-order partial derivative of $p_i$ with respect to a given $\lambda_q$ will bring a term in $k^{n+2}$; this leads to:
$$
\D\forall r \in \N^*: \left|\frac{\partial^{r} p_i}{\partial\lambda_{q_1}\dots\partial\lambda_{q_r}}\right| \leq \left(k^{n+2}\right)^{r}=k^{r(n+2)},
$$
which corresponds to the second inequality of (\ref{deriv_part_ordre_m_pi_vs_lambda_l}).

\end{prooff}
We can now prove the following theorem in order to obtain the estimate for the canonical basis $(p_i)_{i=1,N}$ with respect to the semi-norms $|.|_{l,p,K}$.
\begin{theorem}\label{Estimation_pi}
Let  $\rho$ be the diameter of the largest ball that can be inscribed in $K$. Let $(p_i)_{i=1,N}$ be the canonical basis of $P_k(K)$ defined in (\ref{shape_function}), $(k,l,n)$ three integers and $p$ a positive real number such that:
\begin{equation}\label{parameter_condition_0}
\D k+1 > l + \frac{n}{p},\, (0 < p < +\infty).
\end{equation}
Then, there exists two positive constants $C_0$ and $C_l$ independent of $k$ such that
\begin{equation}\label{Norm_0_2_and_Norm_1_2}
\D\forall p\in\R^*\!\!: |p_i|_{0,p,K} \leq C_0\,k^{n+1} \,\mbox{ and }\,\, \forall l\in\N^{*}\!\!:|p_i|_{l,p,K} \leq C_l\,\frac{k^{l(n+2)}}{\rho^{\,l}}\,.
\end{equation}
\end{theorem}
\begin{prooff} Let us consider the canonical basis of polynomials $(p_i)_{i=1,N}$ of $P_k(K)$ defined by (\ref{shape_function}) and (\ref{P_ij}). \sa
Then, due to remark 2.2 in R. Arcangeli and J. L. Gout \cite{Arcangeli_Gout}, for each polynomial $p_i$, we have for all  $l\geq 0$ for which (\ref{parameter_condition_0}) holds :
\begin{equation}\label{Norm_Grad_L2_pi_0}
\D |p_i|_{l,p,K} \leq \frac{1}{\rho^{\,l}}\left\{\int_{K}\left[\sum_{|\alpha|=l}\frac{l!}{\alpha !}\left|\partial^\alpha p_i(x)\right|\right]^{p}\!\!\!dx\right\}^{\frac{1}{p}} = \frac{1}{\rho^{\,l}}\left\{\int_{K}\left[\sum_{|\alpha|=l}\frac{l!}{\alpha !}\left|\frac{\partial^{|\alpha|} p_i(x)}{\partial x_{1}^{\alpha_1}\dots\partial x_{n}^{\alpha_n}}\right|\right]^{p}\!\!\!dx\right\}^{\frac{1}{p}}\!\!,
\end{equation}
where $\alpha != \alpha_1 !\dots\alpha_n !$ and $\rho$ is the supremum of the diameters of the inscribed spheres within the $n$-simplex $K$.\sa
So, when $l=0$, (\ref{Norm_Grad_L2_pi_0}) together with the first inequality of (\ref{deriv_part_ordre_m_pi_vs_lambda_l}) directly leads to:
\begin{equation}\label{Norm_Grad_L2_pi_0_l=0}
\D |p_i|_{0,p,K} \leq \left\{\int_{K}\left|p_i(x)\right|^{p}dx\right\}^{\frac{1}{p}} \leq \mbox{ mes}(K)^{1/p}\,k^{n+1},
\end{equation}
which corresponds to the first part of (\ref{Norm_0_2_and_Norm_1_2}) with $C_0=\mbox{ mes}(K)^{1/p}$.\sa
Let us now consider the case where $l\geq 1$. Here, each first-order partial derivative $\D \frac{\partial p_i}{\partial x_j}$ can be written as
\begin{equation}\label{D_pi_D_x_j}
\D \frac{\partial p_i}{\partial x_j} = \sum_{q=1}^{n+1}\frac{\partial p_i}{\partial \lambda_q}\frac{\partial \lambda_q}{\partial x_j},
\end{equation}
where $\D \frac{\partial \lambda_q}{\partial x_j}$ is a constant denoted $\Lambda^{q}_j$ which does not depend on $k$, since $\lambda_q$ is a polynomial of degree one, and we rewrite (\ref{D_pi_D_x_j}) as:
\begin{equation}
\D \frac{\partial p_i}{\partial x_j} = \sum_{q=1}^{n+1}\Lambda^{q}_j\frac{\partial p_i}{\partial \lambda_q}.
\end{equation}
Therefore, in the same way, the second-order partial derivatives are given by:
$$
\D \frac{\partial^2 p_i}{\partial x_j\partial x_k} = \sum_{q_1=1}^{n+1}\sum_{q_2=1}^{n+1}\Lambda^{q_1}_{j}\Lambda^{q_2}_{k}\frac{\partial^2 p_i}{\partial\lambda_{q_1}\partial\lambda_{q_2}},
$$
and more generally for any non zero multi-index $\alpha=(\alpha_1,\dots,\alpha_n)$ whose length is denoted $|\alpha|$, we get:
\begin{eqnarray*}
\D \hspace{-0.7cm}\frac{\partial^{|\alpha|} p_i}{\partial x_{1}^{\alpha_1}\dots\partial x_{n}^{\alpha_n}} \hspace{0.3cm} & \, = \, & \,\dots         \nonumber \\[0.2cm]
\D \hspace{-0.6cm}\left(\sum_{q^{1}_1=1}^{n+1}\hspace{-0.2cm}\dots\hspace{-0.2cm}\sum_{q^{1}_{\alpha_1}=1}^{n+1}\hspace{-0.1cm}\right)\hspace{-0.1cm}\dots\hspace{-0.1cm} \left(\sum_{q^{n}_1=1}^{n+1}\hspace{-0.2cm}\dots\hspace{-0.2cm}\sum_{q^{n}_{\alpha_n}=1}^{n+1}\hspace{-0.1cm}\right)
& & \hspace{-1.3cm}\left(\Lambda^{q^{1}_1}_{1}\hspace{-0.1cm}\dots\hspace{-0.1cm}\Lambda^{q^{1}_{\alpha_1}}_{1} \hspace{-0.2cm}\dots\hspace{-0.1cm}\Lambda^{q^{n}_{1}}_{n}\hspace{-0.2cm}\dots\hspace{-0.1cm}\Lambda^{q^{n}_{\alpha_n}}_{n}\right) \hspace{-0.1cm}\frac{\partial^{|\alpha|}p_i}{\left(\partial\lambda_{q^{1}_1}\dots\partial\lambda_{q^{1}_{\alpha_1}}\right) \dots\left(\partial\lambda_{q^{n}_1}\dots\partial\lambda_{q^{n}_{\alpha_n}}\right)}.
\end{eqnarray*}
Now, by using the second estimate of (\ref{deriv_part_ordre_m_pi_vs_lambda_l}) where we set $r=|\alpha|$, this gives the following estimate:
\begin{equation}\label{derive_pi_majo_1}
\D \forall \alpha\in\N^n,\,|\alpha|>0: \left|\frac{\partial^{|\alpha|} p_i(x)}{\partial x_{1}^{\alpha_1}\dots\partial x_{n}^{\alpha_n}}\right| \leq \left[(n+1)\Lambda\right]^{|\alpha|}k^{|\alpha|(n+2)},
\end{equation}
where we set $\D\Lambda\equiv\max_{(j,q)}\Lambda^{q}_{j}$, $(j,q)\in \{1,\dots,n\}\times\{1,\dots n+1\}$ .\sa
Finally, from (\ref{derive_pi_majo_1}) we can derive:
\begin{equation}\label{derive_pi_majo_3}
\D \forall l\in\N^*: \sum_{|\alpha|=l}\frac{l!}{\alpha !}\left|\frac{\partial^{|\alpha|} p_i(x)}{\partial x_{1}^{\alpha_1}\dots\partial x_{n}^{\alpha_n}}\right| \leq \left[(n+1)\Lambda\right]^{l}l!\,n^l k^{l(n+2)},
\end{equation}
since $n^l$ corresponds to the number of partial derivatives of order $l$ in $\R^n$ for the polynomials $p_i$. \sa
Therefore, one can estimate the $|.|_{l,p,K}-$norm for each polynomial $p_i, (1 \leq i\leq N),$ due to (\ref{Norm_Grad_L2_pi_0}) and (\ref{derive_pi_majo_3}), and finally obtain:
\begin{equation}\label{Norm_L2_pi_00_1}
\D \forall l\in\N^*: |p_i|_{l,p,K} \leq \left[\frac{\left[n(n+1)\Lambda\right]^{l} l!\,\mbox{ mes}(K)^{1/p}}{\rho^{\,l}}\right]k^{l(n+2)},
\end{equation}
which corresponds to the second part of (\ref{Norm_0_2_and_Norm_1_2}), with  $C_l=\left[n(n+1)\Lambda\right]^{l} l!\,\mbox{ mes}(K)^{1/p}$.
\end{prooff}
\section{Explicit $k-$dependence in \emph{a priori} $P_k$ finite element error estimates}\label{explicit_estimate}
\vspace{-0.25cm}
\noindent We are now in a position to derive a $k$-explicit dependence of the constant involved in a $W^{m,p}$ \emph{a priori} error estimate for $P_k$ Lagrange finite elements.\\[0.2cm]
This is the purpose of the following theorem:
\begin{theorem}\label{Thm_error_estimate}
Let the hypothesis of C\'ea's Lemma \ref{Lemme_Cea} hold with $W$and $V$ defined by (\ref{V_and_W}). Let $(k,m,n)$ be three integers and $p$ a positive real number satisfying
\begin{eqnarray}
\D \mbox{if }\hspace{0.2cm} \frac{n}{p} < 1 & \mbox{then} & m \leq k, \label{cond_parametre_1}\\[0.2cm]
\D \mbox{if }\hspace{0.2cm} \frac{n}{p} \geq 1 & \mbox{then} & m \leq k-1 \,\mbox{ and }\, k+1-\frac{n}{p}>0. \label{cond_parametre_2}
\end{eqnarray}
Suppose that the approximation $u_h \in W_h$ is a continuous piecewise function composed by polynomials which belong to $P_k(K_{\mu}), (1 \leq \mu \leq N_K)$.
\sa If the exact solution $u$ to problem (\ref{VP}) belongs to $W^{k+1,p}(\Omega)$, the approximation $u_h$, solution to problem (\ref{VP_h}), converges to $u$ in $W^{m,p}(\Omega)$ when $h$ goes to zero, and we have
\vspace{0.1cm}
\begin{equation}\label{estimation_error}
\|u_h-u\|_{m,p,\Omega} \hs \leq \hs \mathscr{C}_k\,h^{k+1-m} \, |u|_{k+1,p,\Omega}\,,
\end{equation}
\vspace{0.1cm}
where $\mathscr{C}_k$ is a positive constant independent of $h$ defined by:
\begin{equation}\label{C_k_Estimation}
\D \mathscr{C}_k \, = \, C
\frac{(k+n)^{n}\,k^{\,m(n+2)}}{(k-m)!\left(\!\!k+\!1-m-\D\frac{n}{p}\right)}
\end{equation}
Above, $C$ is a positive constant which does not depend on $k$.
\end{theorem}
\begin{prooff}
The proof of this theorem is based on the paper of R. Arcangeli and J.L. Gout \cite{Arcangeli_Gout}, itself an extension of the paper by P.G. Ciarlet and P.A. Raviart \cite{Ciarlet_Raviart}. \sa
To this end, let us first recall the conditions of theorem 2.1 of R. Arcangeli and J.L. Gout. \sa
Let $\Omega$ be an open, bounded and non empty convex subset of $\R^n$ and $\Gamma$ its boundary. Let us denote by $P_k$ the space of polynomial functions of degree less than or equal to $k$. We assume that $\Sigma=\{a_i\}_{i=1,N}$ is a $P-$unisolvent set of points which belongs to $\bar{\Omega}$, where $P$ denotes a finite-dimensional space of dimension $N$ composed by functions defined on $\bar{\Omega}$ such that $P_k \subset P \subset C^{k}(\bar{\Omega})$.\sa
Then, for all $u\in W^{k+1,p}(\Omega)$ and for all integer $l \geq 0$ such that
\begin{equation}\label{parameter_condition}
\D k+1 > l + \frac{n}{p},
\end{equation}
we have:
\begin{eqnarray}
\D |u -\Pi_h u|_{l,p,\Omega} & \, \leq \, & \frac{1}{(k-l)!}\frac{1}{\left(k+1-l-\D\frac{n}{p}\right)}|u|_{k+1,p,\Omega}\,h^{k+1-l}  \nonumber \\[0.2cm]
& \, + \, & \frac{1}{(\mbox{mes }\Omega)^{1/p}}\frac{1}{k!} \frac{1}{\left(k+1-\D\frac{n}{p}\right)}\left(\sum_{i=1}^{N}|p_i|_{l,p,\Omega} \right)|u|_{k+1,p,\Omega}\,h^{k+1}, \label{Estimation_Arc_Gout_000}
\end{eqnarray}
where $\Pi_h$ is the classical Lagrange interpolation which consists in interpolating the set of points $\Sigma$ in $\R^n$ by a polynomial function of a given degree $k$, and $(p_i)_{i=1,N}$ are the unique functions such that
$$p_{i}(Mj)=\delta_{ij}, \forall\, M_j \in \Sigma, \forall \, 1\leq i, j\leq N,$$
where $\delta_{ij}$ denotes the classical Kronecker symbol. \sa
First of all, let us remark that, since we are interested in getting an \emph{a priori} error estimate in $W^{m,p}(\Omega)$ for the exact solution $u$ to the variational formulation (\ref{VP}) defined in (\ref{VP}), we will need to write estimates (\ref{Estimation_Arc_Gout_000}) for all values of $l$ between 0 and $m$. It means that condition (\ref{parameter_condition}) also needs to be satisfied from $l=0$ to $m$, which implies that the following inequality must hold true:
$$
\D \frac{n}{p} < k+1-m.
$$
Now, to guarantee this condition, according to the ratio $\D\frac{n}{p}$, we get two conditions, conditions (\ref{cond_parametre_1}) and (\ref{cond_parametre_2}). \sa
Particularly, for the usual case where $p=2$ and $n=2$, condition (\ref{cond_parametre_2}) implies that when considering finite element $P_1$, estimate (\ref{estimation_error}) will only be written for $m=0$ which corresponds to the $L^2$-norm. However, the finite element $P_1$ would also be considered with the $W^{m,p}$-norm by adapting our theorem with another result from R. Arcangeli and J.L. Gout, (see remark 2.3 and theorem 1.1 in \cite{Arcangeli_Gout}). \sa
Thus, for our objectives, we write (\ref{Estimation_Arc_Gout_000}) for the following conditions:
\begin{itemize}
\item $\Omega = K_{\mu}$, $(1\leq \mu \leq N_k)$, a $n$-simplex which belongs to a regular mesh ${\mathcal T}_h$.
\item $u$ is the exact solution in $W^{m,p}(\Omega) \cap W^{k+1,p}(\Omega)$, to the variational formulation (\ref{VP}).
\item The set of points $\Sigma$ in $\R^n$ correspond to the $P_k$ finite element nodes of  $K_{\mu}$.
\item The global interpolate function $\Pi_h u$ is replaced by the local interpolate one $\Pi_{K_{\mu}}u$ on the $n$-simplex $K_{\mu}$.
\end{itemize}
Then, due to (\ref{Norm_L2_pi_00_1}), estimate (\ref{Estimation_Arc_Gout_000}) becomes, $\forall l=1,\dots,m$:
\begin{eqnarray}
\D |u -\Pi_{K_{\mu}} u|_{l,p,K_\mu} & \, \leq \, & \frac{1}{(k-l)!}\frac{1}{\left(k+1-l-\D\frac{n}{p}\right)}|u|_{k+1,p,K_\mu}\,h_{K_{\mu}}^{k+1-l}  \nonumber \\[0.3cm]
& \, + \, & \frac{1}{\rho_{K_\mu}^{l}}\left(\frac{\left[n(n+1)\Lambda\right]^{l}l!}{k!\left(k+1-\D\frac{n}{p}\right)}\frac{(k+n)!}{n!\,k!}\,k^{l(n+2)}\right)
|u|_{k+1,p,K_\mu}\,h_{K_{\mu}}^{k+1}, \label{Estimation_Arc_Gout}
\end{eqnarray}
where we used (\ref{Dim_Pk}). \sa
So, (\ref{Estimation_Arc_Gout}) becomes:
\begin{eqnarray}
\hspace{-0.6cm}\D |u -\Pi_{K_{\mu}} u|_{l,p,K_\mu}\hspace{-0.3cm} & \leq  & \hspace{-0.3cm}\left[\frac{1+\D\left(\frac{\left[n(n+1)\sigma \Lambda\right]^{l}l!\, (k+1)\hspace{-0.1cm}\dots\hspace{-0.1cm}(k+n)\,k^{l(n+2)}}{n!}\right)}{(k-l)!\left(k+1-l-\D\frac{n}{p}\right)}\right] |u|_{k+1,p,K_\mu} h_{K_{\mu}}^{k+1-l}, \nonumber \\[0.3cm]
\hspace{-0.6cm}\hspace{-0.3cm}& \leq & \hspace{-0.3cm}\D\left(\!\!\frac{}{}1+\D\frac{\left[n(n+1)\sigma\right]^{m}\!\Lambda^{\!*}\,m!}{n!}\right) \frac{(k+n)^n\, k^{m(n+2)}}{(k-m)!\,\left(\!\!k+\!1-m-\D\frac{n}{p}\right)}\,\, |u|_{k+1,p,K_\mu} h_{K_{\mu}}^{k+1-l},\label{Estimation_Arc_Gout_m_1}
\end{eqnarray}
where $\D\Lambda^{\!*}\equiv\max_{0\leq l \leq m}\Lambda^l$, and $\sigma$ a given number such that $\sigma\geq 1$ and $\D\frac{h_{K_{\mu}}}{\rho_{K_{\mu}}}\leq \sigma$, $\forall\,K_{\mu} \in {\mathcal T}_h$ which we assumed to be a regular mesh.\sa
For simplicity, we rewrite (\ref{Estimation_Arc_Gout_m_1}) as follows:
\begin{equation}\label{Estimation_Arc_Gout_m_2}
\D |u -\Pi_{K_{\mu}} u|_{l,p,K_\mu} \, \leq \, \D C_1(\sigma,\Lambda^{\!*},m,n)\, \frac{(k+n)^n\, k^{m(n+2)}}{(k-m)!\,\left(\!\!k+\!1-m-\D\frac{n}{p}\right)}\, |u|_{k+1,p,K_\mu}\,h_{K_{\mu}}^{k+1-l},
\end{equation}
where we introduced constant $C_1(\sigma,\Lambda^{\!*},m,n)$ defined by:
\begin{equation}\label{xi}
\D C_1(\sigma,\Lambda^{\!*},m,n) \, \equiv \, 1+\D\frac{\left[n(n+1)\sigma\right]^{m}\!\Lambda^{\!*}\,m!}{n!}\,.
\end{equation}
Now, when $l=0$, due to (\ref{Norm_Grad_L2_pi_0_l=0}), estimate (\ref{Estimation_Arc_Gout_000}) becomes:
\begin{eqnarray*}
\D |u -\Pi_{K_{\mu}} u|_{0,p,K_{\mu}} & \, \leq \, & \frac{1}{k!}\frac{1}{\left(k+1-\D\frac{n}{p}\right)}|u|_{k+1,p,K_{\mu}}\,h^{k+1}  \nonumber \\[0.2cm]
& \, + \, & \frac{1}{k!} \frac{1}{\left(k+1-\D\frac{n}{p}\right)}\frac{(k+n)!}{n!\,k!}\,k^{n+1}|u|_{k+1,p,K_{\mu}}\,h^{k+1},
\end{eqnarray*}
which leads to:
\begin{equation}\label{base_caracteristique}
\D |u -\Pi_{K_{\mu}}|_{0,p,K_{\mu}} \leq C_2(n) \frac{(k+n)^nk^{m(n+2)}}{(k-m)!\left(k+1-m-\D\frac{n}{p}\right)}|u|_{k+1,p,K_{\mu}}\,h^{k+1},
\end{equation}
for all $k \geq 1$ and $m \geq 1$, and where we introduced constant $C_2(n)$ defined by:
\begin{equation}\label{C2_n}
\D C_2(n) = 1+\frac{1}{n!}
\end{equation}
Therefore, by the help of (\ref{Estimation_Arc_Gout_m_2})-(\ref{xi}) and (\ref{base_caracteristique})-(\ref{C2_n}), we get the following $W^{m,p}$ local interpolation error estimate:
\begin{eqnarray}
\hspace{-0.8cm}\D \|u -\Pi_{K_{\mu}} u\|^p_{m,p,K_{\mu}} \hspace{-0.2cm} & = & \hspace{-0.1cm}\sum_{l=0}^{m}|u -\Pi_{K_{\mu}}|^{p}_{l,p,K_{\mu}}, \nonumber \\
\hspace{-0.4cm} & \leq & \hspace{-0.1cm} \sum_{l=0}^{m}C^p(\sigma,\Lambda^*,m,n)\hspace{-0.15cm}\left[\frac{(k+n)^n\, k^{m(n+2)}}{(k-m)!\,\left(\!\!k+\!1-m-\D\frac{n}{p}\right)}\!\right]^p
\hspace{-0.3cm}|u|^{p}_{k+1,p,K_\mu} h_{K_{\mu}}^{p(k+1-l)},  \label{Bid_0}
\end{eqnarray}
where constant $C(\sigma,\Lambda^*,m,n)$ is defined by : $\D C(\sigma,\Lambda^*,m,n) =  \max\left(\frac{}{}\!\!C_1(\sigma,\Lambda^{\!*},m,n), C_2(n)\right)$.
Then, (\ref{Bid_0}) leads to:
\begin{equation}\label{Lambda_0000}
\D \|u -\Pi_{K_{\mu}} u\|^p_{m,p,K_{\mu}} \leq \hspace{-0.1cm} C^p(\sigma,\Lambda^*,m,n,p,h)\hspace{-0.15cm}\left[\frac{(k+n)^n\, k^{m(n+2)}}{(k-m)!\,\left(\!\!k+\!1-m-\D\frac{n}{p}\right)}\!\right]^p \hspace{-0.3cm} |u|^{p}_{k+1,p,K_\mu} h^{p(k+1-m)},
\end{equation}
where $\D h \equiv \max_{K_{\mu} \in {\mathcal T}_h} h_{K_{\mu}}$ and $C(\sigma,\Lambda^*,m,n,p,h)\equiv \xi(m,p,h).C(\sigma,\Lambda^*,m,n)$ with $\xi(m,p,h)$ defined as follows:
\begin{equation}\label{varphi}
\D \xi(m,p,h) \equiv \left |
\begin{array}{ll}
\D\hs \left[\frac{1-h^{p(m+1)}}{1-h^p}\right]^{\frac{1}{p}} & \mbox{ if } \hs h\neq 1, \medskip \\
\hs (m+1)^{\frac{1}{p}} & \mbox{ if } \,h=1.
\end{array}
\right.
\end{equation}
Since the mesh ${\mathcal T}_h$ is regular, by the help of (\ref{Lambda_0000}), we get for the whole domain $\Omega$ the following global interpolation error estimate:
\begin{eqnarray}
\D \|u -\Pi_h u\|_{m,p,\Omega} & = & \hspace{-0.2cm} \left(\sum_{K_{\mu}\in {\mathcal T}_h}\|u -\Pi_{K_{\mu}}u\|_{m,p,K_{\mu}}^p\!\!\right)^{\!\!1\!/p} \nonumber \\
& \leq & \hspace{-0.2cm} C(\sigma,\Lambda^*,m,n,p,h)\left[\frac{(k+n)^n k^{m(n+2)}}{(k-m)!\,\left(\!\!k+\!1-m-\D\frac{n}{p}\right)}\!\right]\left(\!\sum_{K_{\mu}\in {\mathcal T}_h}|u|^p_{k+1,p,K_{\mu}}\!\!\right)^{\!\!1\!/p}\!\!\!\!\!\!h^{k+1-m},\nonumber\\[0.2cm]
& \leq & \hspace{-0.2cm} C(\sigma,\Lambda^*,m,n,p,h)\left[\frac{(k+n)^n k^{m(n+2)}}{(k-m)!\,\left(\!\!k+\!1-m-\D\frac{n}{p}\right)}\!\right]|u|_{k+1,p,\Omega}\,h^{k+1-m}. \label{Error_estimate_for_any_k}
\end{eqnarray}
Then, estimate (\ref{C_k_Estimation}) is proved, provided that one takes into account estimate (\ref{Cea_Banach}) of C\'ea's Lemma \ref{Lemme_Cea}. Indeed, consider the $W^{m,p}-$norm to measure the difference between the exact solution $u$ to the variational problem (\ref{VP}) and its approximation $u_h$ solution to (\ref{VP_h}), we have:
\begin{equation}\label{Erreur_approx_Erreur_global_interpol}
\|u-u_h\|_{m,p,\Omega} \hs \leq \hs \left(1+ \frac{\|a\|_{W^{m,p},W^{m'\!,p\,'}}}{\alpha_h}\right)\,\|u-\Pi_h u\|_{m,p,\Omega},
\end{equation}
where we choose in (\ref{Cea_Banach}) $w_h\in W_h$ equal to $\Pi_h u$, $\|a\|_{W^{m,p},W^{m'\!,p\,'}}$ as defined in (\textbf{BNB}), and $\alpha_h$ being the constant of the discrete inf-sup condition, see  (\textbf{BNB1$_h$}). \sa
Then, replacing expression (\ref{Error_estimate_for_any_k}) in inequality (\ref{Erreur_approx_Erreur_global_interpol}) leads to:
\begin{equation}\label{Estim_Preuve}
\|u-u_h\|_{m,p,\Omega} \hs \leq \hs \left(1+ \frac{\|a\|_{W^{m,p},W^{m'\!,p\,'}}}{\alpha_h}\right)\,C(\sigma,\Lambda^*,m,n,p,h)\left[\frac{(k+n)^n k^{m(n+2)}}{(k-m)!\,\left(\!\!k+\!1-m-\D\frac{n}{p}\right)}\!\right]|u|_{k+1,p,\Omega}\,h^{k+1-m}.
\end{equation}
Now, since $\xi(m,p,h)$ introduced in (\ref{varphi}) is bounded as $h\leq diam\,({\bar{\Omega}})$, $C(\sigma,\Lambda^*,m,n,p,h)$ is uniformly bounded with $h$. Hence, there exists $C(\sigma,\Lambda^*,m,n,p)$ independent of $h$ such that:
$$
C(\sigma,\Lambda^*,m,n,p,h) \leq C(\sigma,\Lambda^*,m,n,p).
$$
Consequently, by defining constant $\mathscr{C}_k$ by
\begin{equation}\label{C_k_Value}
\D \mathscr{C}_k \equiv \left(1+ \frac{\|a\|_{W^{m,p},W^{m'\!,p\,'}}}{\alpha_h}\right)C(\sigma,\Lambda^*,m,n,p)\,\frac{(k+n)^n k^{m(n+2)}}{(k-m)!\,\left(\!\!k+\!1-m-\D\frac{n}{p}\right)},
\end{equation}
we obtain the error estimate (\ref{estimation_error})-(\ref{C_k_Estimation}), with $C= \left(1+ \D\frac{\|a\|_{W^{m,p},W^{m'\!,p\,'}}}{\alpha_h}\right)C(\sigma,\Lambda^*,m,n,p)$.
\end{prooff}
\section{Application to relative finite elements accuracy}\label{FEM_accuracy}
\noindent In this section, we apply inequality (\ref{estimation_error}) of Theorem \ref{Thm_error_estimate} to evaluate the relative accuracy between two finite elements.
Hereafter, we will replace notation $u_h$ with $u_h^{(k)}$, in order to highlight the degree $k$ of the polynomials involved in $P_{k}(K_{\mu})$.\sa
In \cite{ChasAs18}, regarding a problem set in the usual Sobolev space $H^1(\Omega)$, we introduced a probabilistic framework which enables one to compare the relative accuracy of two finite elements of different degrees in a non standard way. Indeed, we claimed that quantitative uncertainties exist in the approximate solution $u^{(k)}_h$, due for instance to the quantitative uncertainties that are commonly produced in the mesh generation.\sa
For this reason, we have considered the approximation error as a random variable, and we aimed at evaluating the probability of the difference between two $H^{1}-$approximation errors of $u- u_h^{(k_1)}$ and $u- u_h^{(k_2)}$ corresponding to finite elements $P_{k_1}$ and $P_{k_2}, (k_1<k_2)$.\sa
Here, in the same way, one can only infer that the value of the approximation error $\|u^{(k)}_h-u\|_{m,p,\Omega}$ belongs to the interval $[0,\mathscr{C}_k |u|_{k+1,p,\Omega}\,h^{k+1-m}]$, using error estimates (\ref{estimation_error})-(\ref{C_k_Estimation}).\sa
As a consequence, for fixed values of $k, m$ and $p$, we define the following random variable $X_{m,p}^{(k)}$ by:
\begin{eqnarray*}
X_{m,p}^{(k)} : & {\bf\Omega} & \hspace{0.1cm}\rightarrow \hspace{0.2cm}[0,\mathscr{C}_k |u|_{k+1,p,\Omega}\,h^{k+1-m}] \noindent \\
& \boldsymbol{\omega}\equiv u^{(k)}_h & \hspace{0.1cm} \mapsto \hspace{0.2cm}\D X_{m,p}^{(k)}(\boldsymbol{\omega}) = X_{m,p}^{(k)}(u^{(k)}_h) = \|u^{(k)}_h-u\|_{m,p,\Omega},
\end{eqnarray*}
where the probability space ${\bf\Omega}$ contains all the possible results for a given random trial, namely, all possible grids that the involved meshing tool can generate for a given value of $h$. Equivalently, ${\bf\Omega}$ consists of all the possible corresponding approximations $u^{(k)}_h$. Below, for simplicity, we will set: $X_{m,p}^{(k)}(u^{(k)}_h)\equiv X_{m,p}^{(k)}(h)$. \sa
Now, regarding the absence of information concerning the more likely or less likely values of norm $\|u^{(k)}_h-u\|_{m,p,\Omega}$ within the interval $[0, \mathscr{C}_k |u|_{k+1,p,\Omega}\,h^{k+1-m}]$, we assume that the random variable $X_{m,p}^{(k)}(h)$ has a uniform distribution on the interval $[0, \mathscr{C}_k |u|_{k+1,p,\Omega}\,h^{k+1-m}]$ in the following sense:
$$
\forall (\alpha,\beta), 0 \leq \alpha < \beta \leq \mathscr{C}_k |u|_{k+1,p,\Omega}\,h^{k+1-m}: Prob\left\{X_{m,p}^{(k)}(h) \in [\alpha,\beta]\right\}=\frac{\beta-\alpha}{\mathscr{C}_k |u|_{k+1,p,\Omega}\,h^{k+1-m}}.
$$
The above equation means that if one slides interval $[\alpha,\beta]$ anywhere in $[0, \mathscr{C}_k |u|_{k+1,p,\Omega}\,h^{k+1-m}]$, the probability of the event $\D\left\{X_{m,p}^{(k)}(h) \in [\alpha,\beta]\right\}$ does not depend on the localization of  $[\alpha,\beta]$ in $[0, \mathscr{C}_k |u|_{k+1,p,\Omega}\,h^{k+1-m}]$, but only on its length; this reflects the property of  uniformity for $X_{m,p}^{(k)}$. \sa
Hence,  it is straightforward to extend the theorem proved in \cite{ChasAs18} for the $H^{1}$ case to the $W^{m,p}$ context. This yields the following result, which estimates the probability of event $\D\left\{X_{m,p}^{(k_2)}(h) \leq X_{m,p}^{(k_1)}(h)\right\}$. \sa
Let $C_{k_{i}}$ be equal to $C_{k_i} \equiv \mathscr{C}_{k_i} |u|_{k_{i}+1,p,\Omega}$, for $i=1,2$, and let  $h^*_{m,p}$ be defined as:
\begin{equation}\label{h*}
\D h^*_{m,p} \equiv \left( \frac{C_{k_{1}}}{C_{k_{2}}} \right)^{\frac{1}{k_2-k_1}}.
\end{equation}
As in \cite{ChasAs18}, by changing the $H^1-$norm to the $W^{m,p}$ one, we can derive that:
\begin{equation}\label{Nonlinear_Prob}
\D Prob\left\{ X_{m,p}^{(k_2)}(h) \leq X_{m,p}^{(k_1)}(h)\right\} = \left |
\begin{array}{ll}
\D \hs 1 - \frac{1}{2}\!\left(\!\frac{\!\!h}{h^*_{m,p}}\!\right)^{\!\!k_2-k_1} & \mbox{ if } \hs 0 < h \leq h^*_{m,p}, \medskip \\
\D \hs \frac{1}{2}\!\left(\!\frac{h^*_{m,p}}{\!\!h}\!\right)^{\!\!k_2-k_1} & \mbox{ if } \hs h \geq h^*_{m,p}.
\end{array}
\right.
\end{equation}

\noindent Then, using (\ref{C_k_Value}), one can rewrite $h^*_{m,p}$ defined in (\ref{h*}) as follows:
\begin{equation}\label{explicit_h*}
\D h^*_{m,p} = \left[
\frac{\left(1+ \D\frac{\|a\|_{W^{m,p},W^{m'\!,p\,'}}}{\alpha_{h,k_1}}\right)}{\left(1+ \D\frac{\|a\|_{W^{m,p},W^{m'\!,p\,'}}}{\alpha_{h,k_2}}\right)}
\left(\frac{k_1+n}{k_2+n}\right)^{n}\left(\frac{k_1}{k_2}\right)^{m(n+2)}\frac{(k_2-m)!}{(k_1-m)!}\,\frac{\left(k_2+1-m-\D\frac{n}{p}\right)}{\left(k_1+1-m-\D\frac{n}{p}\right)}\,\frac{|u|_{k_{1}+1,p,\Omega} }{|u|_{k_{2}+1,p,\Omega}}
\right]^{\frac{1}{k_2-k_1}}\,,
\end{equation}
where $\alpha_{h,k_1}$ and $\alpha_{h,k_2}$ denotes the $\alpha_h$ appearing in generalized C\'ea's Lemma \ref{Lemme_Cea}, associated to finite elements $P_{k_1}$ and $P_{k_2}$, respectively.
\begin{remark}
Notice that, as proposed in \cite{Arxiv2}, one can derive another law of probability to evaluate the most accurate finite element between $P_{k_{1}}$ and $P_{k_{2}}$. More precisely, for $h<h^*_{m,p}$, assuming the independence of events $A\equiv\D\left\{X_{m,p}^{(k_2)}(h) \leq X_{m,p}^{(k_1)}(h)\right\}$ and $B\equiv \left\{X_{m,p}^{(k_1)}(h) \in [C_{k_2} h^{k_2},C_{k_1} h^{k_1}]\right\}$, one can obtain the following law of probability:
\begin{equation}\label{Heaviside_Prob}
\D Prob\left\{ X_{m,p}^{(k_2)}(h) \leq X_{m,p}^{(k_1)}(h)\right\} = \left |
\begin{array}{ll}
\hs 1 & \mbox{ if } \hs 0 < h < h^*_{m,p}, \medskip \\
\hs 0 & \mbox{ if } \hs h> h^*_{m,p}\,.
\end{array}
\right.
\end{equation}
\end{remark}
The probability distribution (\ref{Heaviside_Prob}) is obtained by replacing the uniform distribution assumption in (\ref{Nonlinear_Prob}) by the independence of events $A$ and $B$. However, with no prior information about the independence of these events, the more "natural" probabilistic law is (\ref{Nonlinear_Prob}).\sa
Therefore, in what follows, we take a fixed value for $k_1$ (that we will denote $k$ in the sequel),  and we study the asymptotic behavior of the accuracy between $P_k$ and $P_{k+q}$, when $q$ goes to $+\infty$: this will give us the asymptotic relation between the two probabilistic laws (\ref{Nonlinear_Prob}) and (\ref{Heaviside_Prob}). \sa
To this end, it is convenient to introduce notation $\D\left(\mathcal{P}_{q}(h)\right)_{q \in \N^\star}$  corresponding to the sequence of functions defined by (\ref{Nonlinear_Prob}), namely:
\begin{equation}\label{P(h)}
\D\forall q\in\N^*:\mathcal{P}_{q}(h) \equiv Prob\left\{ X_{m,p}^{(k+q)}(h) \leq X_{m,p}^{(k)}(h)\right\} = \left |
\begin{array}{ll}
\D \hs 1 - \frac{1}{2}\!\left(\!\frac{\!\!h}{h^{*}_q}\!\right)^{\!q} & \mbox{ if } \hs 0 < h \leq h^{*}_q, \medskip \\
\D \hs \frac{1}{2}\!\left(\!\frac{h^{*}_q}{\!\!h}\!\right)^{\!q} & \mbox{ if } \hs h \geq h^{*}_q\,.
\end{array}
\right.
\end{equation}
Above, we denote by $h^{*}_q$ the $h^*_{m,p}$ expressed as a function of $q$ for given values of $k,m$ and $p$, that is:
\begin{equation}\label{h*q}
\hspace{-1cm}\D h^*_q = \left[
\frac{\left(1+ \D\frac{\|a\|_{W^{m,p},W^{m'\!,p\,'}}}{\alpha_{h,k}}\right)}{\left(1+ \D\frac{\|a\|_{W^{m,p},W^{m'\!,p\,'}}}{\alpha_{h,k+q}}\right)}
\left(\frac{k+n}{k+q+n}\right)^{n}\left(\frac{k}{k+q}\right)^{m(n+2)}\frac{(k+q-m)!}{(k-m)!}\,\frac{\left(k+q+1-m-\D\frac{n}{p}\right)}{\left(k+1-m-\D\frac{n}{p}\right)}
\,\frac{|u|_{k+1, p,\Omega} }{|u|_{k+q+1, p,\Omega}}\right]^{\frac{1}{q}}\hspace{-0.1cm}.
\end{equation}
To obtain the asymptotic behavior of sequence $\D\left(\mathcal{P}_{q}(h)\right)_{q \in \N^\star}$, we first have to compute the limit of sequence $\D\left(h^*_q\right)_{q\in\N}$. \sa
It is the purpose of the following lemma:
\begin{lemma}\label{Conv_asympt_h*q}
Let $u\in W^{r,p}(\Omega), (\forall\, r\in \N),$ be the solution to problem (\ref{VP}) and $(h^*_{q})_{q\in\N^\star}$ the sequence defined by (\ref{h*q}). We assume that sequence $(\alpha_{h,k+q})_{q\in\N^\star}$ satisfies:
\begin{equation}\label{alpha_asympt}
\D\forall  k\in\N, \lim_{q\rightarrow +\infty}\alpha_{h,k+q} = \alpha_{h,k}^* \in\R^*.
\end{equation}
Let $k, m$ and $p$ be fixed such that (\ref{cond_parametre_1}) or (\ref{cond_parametre_2}) holds. \sa
If
\begin{equation}\label{Cond_ Ratio_Semi_Norm}
\D \lim_{q\rightarrow +\infty}\frac{|u|_{k+q+2,p,\Omega}}{|u|_{k+q+1,p,\Omega}} = l, (l\in\R^{*}_{+}),
\end{equation}
then,
\begin{equation}\label{lim_hq*}
\lim_{q\rightarrow +\infty}h^*_q = +\infty.
\end{equation}
\end{lemma}
\noindent \begin{prooff}
From (\ref{h*q}), we readily get:
\begin{equation}\label{C_k_0}
\hspace{-1cm}\D\left(h^*_q\right)^q
=
\frac{\left(1+ \D\frac{\|a\|_{W^{m,p},W^{m'\!,p\,'}}}{\alpha_{h,k}}\right)}{\left(1+ \D\frac{\|a\|_{W^{m,p},W^{m'\!,p\,'}}}{\alpha_{h,k+q}}\right)}
 \frac{(k+n)^n k^{m(n+2)}}{(k-m)!\,\left(k+1-m-\D\frac{n}{p}\right)} \, \frac{(k+q-m)!\left(k+q+1-m-\D\frac{n}{p}\right)}{\D\left(k+q+n\right)^n \left(k+q\right)^{m(n+2)}}\,.\frac{|u|_{k+1,p,\Omega}}{|u|_{k+q+1,p,\Omega}}.
\end{equation}
Let us first remark that condition (\ref{alpha_asympt}) implies that the following ratio, based on the constant involved in (\ref{Cea_Banach}) of Lemma \ref{Lemme_Cea}, is uniformly bounded and stays strictly positive for any value of $q$. In particular, we have:
\begin{equation}\label{rapp_alpha}
\D \lim_{q\rightarrow +\infty}\frac{\left(1+ \D\frac{\|a\|_{W^{m,p},W^{m'\!,p\,'}}}{\alpha_{h,k}}\right)}{\left(1+ \D\frac{\|a\|_{W^{m,p},W^{m'\!,p\,'}}}{\alpha_{h,k+q}}\right)} = \beta_{h,k}^* \in\R^{*}.
\end{equation}
Then, using Stirling's formula when $q$ goes to $+\infty$, we first remark that
\begin{eqnarray}
\frac{(k\!+\!q\!-\!m)!\left(k\!+\!q\!+\!1\!-\!m\!-\!\D\frac{n}{p}\right)}{\D\left(k+q\right)^{m(n+2)}\left(k+q+n\right)^n} \hspace{-0.3cm}& \underset{q \rightarrow +\infty}{\sim} & \hspace{-0.3cm}\frac{\sqrt{2\pi (k+q-m)}\D\left(\frac{k+q-m}{e}\right)^{\!\!(k+q-m)}\!\!\!\left(k\!+\!q\!+\!1\!-\!m\!-\!\D\frac{n}{p}\right)}{\left(k+q\right)^{m(n+2)}\left(\!\!\frac{}{}k+q+n\right)^n}, \nonumber \\[0.2cm]
\hspace{-0.3cm}& \underset{q \rightarrow +\infty}{\sim} & \hspace{-0.3cm}\frac{\sqrt{2 \pi}(k+q-m)^{(k+q-m+\frac{1}{2})}}{e^{k+q-m}}\frac{1}{\left(\!\!\frac{}{}k+q\right)^{n+m(n+2)-1}}, \nonumber \\[0.2cm]
\hspace{-0.3cm}& \underset{q \rightarrow +\infty}{\sim} & \hspace{-0.3cm}\sqrt{2 \pi}\,\frac{(k+q)^{k+q-3m-n(m+1)+\frac{3}{2}}}{e^{k+q}}, \label{h*q_Asympt}
\end{eqnarray}
where, according to Euler's formula \cite{Euler}, we have used the following equivalence
$$
\D (k+q-m)^{(k+q-m+\frac{1}{2})} \underset{q \rightarrow +\infty}{\sim} e^{-m}(k+q)^{(k+q-m+\frac{1}{2})}.
$$
Then, substituting (\ref{h*q_Asympt}) in (\ref{C_k_0}) allows us to determine equivalent of $h^{*}_q$ when $q \rightarrow +\infty$.  Using (\ref{rapp_alpha}), one obtains
\begin{equation}\label{Truc}
\D\left(h^*_q\right)^q \underset{q \rightarrow +\infty}{\sim} \Theta \,\frac{\left(1+ \D\frac{\|a\|_{W^{m,p},W^{m'\!,p\,'}}}{\alpha_{h,k}}\right)}{\left(1+ \D\frac{\|a\|_{W^{m,p},W^{m'\!,p\,'}}}{\alpha_{h,k+q}}\right)} e^{-(k+q)}(k+q)^{k+q-3m-n(m+1)+\frac{3}{2}}\,.\frac{|u|_{k+1,p,\Omega}}{|u|_{k+q+1,p,\Omega}},
\end{equation}
where $\Theta$ denotes a constant independent of $q$ defined  by
$$
\Theta \equiv \sqrt{2 \pi}\frac{(k+n)^n k^{m(n+2)}}{(k-m)!\,\left(k+1-m-\D\frac{n}{p}\right)}.
$$
We now introduce two sequences $(v_q)_{q\in\N}$ and $(w_q)_{q\in\N}$, as follows:
$$
\forall\, q\in\N: v_q\equiv \ln\,|u|_{k+q+1,p,\Omega}, \hs w_q \equiv q\,.
$$
Then, owing to condition (\ref{Cond_ Ratio_Semi_Norm}), if sequence $r_q$ defined by the ratio
$$
\D r_q \equiv \frac{v_{q+1}-v_q}{w_{q+1}-w_q} = \ln\left(\frac{|u|_{k+q+2,p,\Omega}}{|u|_{k+q+1,p,\Omega}}\right)
$$
has a limit $L\equiv \ln\, l \in\R$, when $q$ goes to $+\infty$, then $\D\lim_{q\rightarrow +\infty}r_q = L$.\sa
As a consequence, due to the Stolz-Cesaro theorem, (see \cite{OviFur}, p.263-266), the ratio $\D\frac{v_q}{w_q}$ also converges to the same limit $L$ when $q$ goes to $+\infty$:
$$
\D \lim_{q\rightarrow +\infty}\D\frac{v_q}{w_q} = \lim_{q\rightarrow +\infty}\D\frac{\ln\,|u|_{k+q+1,p,\Omega}}{q} = L,
$$
and
\begin{equation}\label{exp(-L)}
\D \lim_{q\rightarrow +\infty}\D\left(\frac{|u|_{k+1,p,\Omega}}{|u|_{k+q+1,p,\Omega}}\right)^{\frac{1}{q}} = \lim_{q\rightarrow +\infty}\D\left(\frac{1}{|u|_{k+q+1,p, \Omega}}\right)^{\frac{1}{q}} =e^{-L}=\frac{1}{l}.
\end{equation}
As a result, from (\ref{rapp_alpha}), (\ref{Truc}) and (\ref{exp(-L)}), one can conclude that $\D h^*_q \underset{+\infty}{\sim} \frac{1}{e\,l}\,q$, which proves (\ref{lim_hq*}).
\end{prooff}
%
%
\begin{remark}$\frac{}{}$ Let us comment on the assumptions of this lemma.
\begin{enumerate}
\item  The hypothesis on the ratio of norms in Eq. (\ref{Cond_ Ratio_Semi_Norm}) might appear very {\em ad-hoc}. Nevertheless, one can easily check, based on several examples, that it is satisfied. Take for instance $u$, solution to a standard Laplace problem solved in a regular domain $\Omega \in \R^2$ (for example a square), with a given regular Dirichlet boundary condition on the boundary $\partial \Omega$ and a regular enough right-hand side, (for details see \cite{Arxiv2}).
\item As a matter of fact, inequality (\ref{rapp_alpha}) is fulfilled in most cases. For instance, assuming the bilinear form $a$ is coercive, the functional framework is necessarily Hilbertian (see Remark 2.3 in \cite{Ern_Guermond}). Since we have considered that $W\equiv W^{m,p}(\Omega) \mbox{ and } V\equiv W^{m',p'}(\Omega)$, with $\frac{1}{p}+\frac{1}{p'}=1$, then $p=2$ and $W=V=H^{m}(\Omega)$.\sa
 So, if we denote by $\alpha$ the coercivity constant and by $\|a\|$ the continuity constant, inequality (\ref{Cea_Banach}) of C\'ea's Lemma \ref{Lemme_Cea} can be expressed with constant $\frac{\|a\|}{\alpha}$ instead of $\left(1+\frac{\|a\|_{W,V}}{\alpha_h}\right)$ and the ratio in the limit (\ref{rapp_alpha}) equals 1 and $\beta_{h,k}^*$ too.
\item In terms of linear algebra, i.e. considering the matrix $\A$ associated with the bilinear form $a$, it can be shown that $\alpha_{h,k}$ (or $\alpha_{h,k+q}$) of (\ref{rapp_alpha}) is related to the smallest eigenvalue of the square matrix $\A^t \A$, i.e. the smallest singular value of $\A$, (see \cite{Ern_Guermond}, Remark 2.23, (iii)). Hence, inequality (\ref{rapp_alpha}) could be checked if one is able to get information about the singular value decomposition of $\A$.
\end{enumerate}
\end{remark}
\noindent We now consider the convergence of sequence $\D\left(\mathcal{P}_{q}(h)\right)_{q \in \N^\star}$ as $q\rightarrow +\infty$. As we will see, due to the definition (\ref{P(h)}) of sequence $\D\left(\mathcal{P}_{q}(h)\right)_{q \in \N^\star}$, pointwise convergence presents a discontinuity at point $h=h^*_q$. Indeed, when $q$ goes to $+\infty$, thanks to lemma \ref{Conv_asympt_h*q}, $h^*_q$ also goes to $+\infty$, and this discontinuity is therefore at $+\infty$. \sa
Thus, to handle this singular behavior, we introduce the weak convergence of the sequence $\D\left(\mathcal{P}_{q}(h)\right)_{q \in \N^\star}$, {\it i.e.} convergence on the sense of distributions. \sa
For the sake of exhaustivity, we briefly recall here some basic notions about distribution theory \cite{Schwartz}, that allows us in passing to introduce the notations we will use. A well-informed reader may skip these few lines.\sa
We denote by $\mathcal{D}(\R)$ the space of functions $\,C^{\infty}(\R)$ with a compact support in $\R$, and by $\mathcal{D'}(\R)$ the space of distributions defined on $\R$. As we will carry out our analysis for all x $\in \R$, we extend the sequence of functions $\D\left(\mathcal{P}_{q}(h)\right)_{q \in \N^\star}$ on $]-\infty, 0[$ by setting: $\forall h \leq 0:\mathcal{P}_{q}(h)=0$.\sa
Therefore, the sequence of extended functions $\D\left(\frac{}{}\!\!\mathcal{P}_{q}(h)\right)_{q\in\N^\star}$ belongs to the space $L^{1}_{loc}(\R)$.\footnote{the space of functions locally integrable for any compact $K$ of $\R$.}
Hence, $\forall\, q\in\N^*$, each function $\frac{}{}\!\!\mathcal{P}_{q}(h)$ can be associated to its regular distribution $T_{\mathcal{P}_{q}}$ defined by:
\begin{equation}\label{P(h)_Distrib}
\D\forall \varphi \in {\cal D(\R)} : \hs <T_{\mathcal{P}_{q}}, \varphi > \hs \equiv \hs \int_{\R}\mathcal{P}_{q}(h)\varphi(h)dh.
\end{equation}
For what follows, we will also need the Heaviside distribution $T_H$ defined by:
$$
\D\forall \varphi \in {\cal D(\R)} : \hs <T_H, \varphi > \hs \equiv \hs \int_{\R}H(h)\varphi(h)dh \hs = \hs \int_{0}^{+\infty}\!\!\varphi(h)dh, $$
where $H(h)=1$ if $h > 0$,  and zero otherwise.
We are now in a position to state  the convergence result of the sequence of distributions $\left(T_{\mathcal{P}_{q}}\right)_{q\in\N^*}$
 in $D'(\R)$.
\begin{theorem}\label{Conv_Simple}
With the same assumptions on $u$ as in Lemma \ref{Conv_asympt_h*q}, let $\D\left(T_{\mathcal{P}_{q}}\right)_{q\in\N^*}$ be the sequence of distributions defined by (\ref{P(h)_Distrib}) and (\ref{P(h)})-(\ref{h*q}).
Then, $\D\left(T_{\mathcal{P}_{q}}\right)_{q\in\N^*}$ converges with respect to the weak-* topology on $D'(\R)$ to the Heaviside distribution $T_H$.
\end{theorem}
\noindent \begin{prooff}$\frac{}{}$
By definition \cite{Schwartz} of the weak convergence in $D'(\R)$, we have to evaluate the limit of the numerical sequence $\D\left(<T_{\mathcal{P}_{q}}, \varphi >\right)_{q\in\N^*}$ when $q$ goes to $+\infty$.\sa
Hence, due to (\ref{P(h)_Distrib}) and (\ref{P(h)}), we have, $\D\forall \varphi \in {\cal D(\R)}:$
\begin{eqnarray}\label{limit_distrib}
\hspace{-0.6cm}\D <T_{\mathcal{P}_{q}}, \varphi > \hspace{-0.1cm}& \equiv & \hspace{-0.1cm}\int_{\R}\mathcal{P}_{q}(h)\varphi(h)dh = \int_{0}^{h^{*}_q}\left[1 - \frac{1}{2}\!\left(\!\frac{\!\!h}{h^{*}_q}\!\right)^{\!q}\right]\varphi(h)dh + \int_{h^{*}_q}^{+\infty}\frac{1}{2}\!\left(\!\frac{h^{*}_q}{\!\!h}\!\right)^{\!q}\varphi(h)dh \nonumber \label{}\\[0.2cm]
\hspace{-0.6cm} & = & \hspace{-0.1cm}\int_{-\infty}^{+\infty}\left[1 - \frac{1}{2}\!\left(\!\frac{\!\!h}{h^{*}_q}\!\right)^{\!q}\right]1\hspace{-0.15cm}1_{[0,h^*_q]}(h)\varphi(h)dh + \int_{-\infty}^{+\infty}\frac{1}{2}\!\left(\!\frac{h^{*}_q}{\!\!h}\!\right)^{\!q}1\hspace{-0.15cm}1_{[h^*_q,+\infty[}(h)\varphi(h)dh,\label{limit_distrib}
\end{eqnarray}
where $1\hspace{-0.15cm}1_{[a,b]}$ denotes the indicator function of interval $[a,b], \forall (a,b)\in \R^2$.\sa
Therefore, to compute the limit of $<T_{\mathcal{P}_{q}}, \varphi >$ when $q$ goes to $+\infty$, we will check the hypothesis of the dominated convergence theorem \cite{Brezis} for the integrals involved in (\ref{limit_distrib}).\sa
For the first one, introduce the sequence of functions $\D\left(\psi_{q}\right)_{q\in\N^*}$ defined on $\R$ by:
$$
\D\forall \, q \in \N^*: \psi_q(h)= \left[1 - \frac{1}{2}\!\left(\!\frac{\!\!h}{h^{*}_q}\!\right)^{\!q}\right]1\hspace{-0.15cm}1_{[0,h^*_q]}(h)\varphi(h).
$$
Then, sequence $\D\left(\psi_{q}\right)_{q\in\N^*}$ exhibits the following properties:
\begin{itemize}
\item It converges pointwise on $\R$ to function $H\varphi$, thanks to the following properties:
\begin{eqnarray*}
\D \forall h \in\R & : &\lim_{q\rightarrow +\infty}1\hspace{-0.15cm}1_{[0,h^*_q]}(h)=1\hspace{-0.15cm}1_{[0,+\infty[}(h), \\[0.2cm]
\D\forall h,\, 0 < h < h^*_q & : & \lim_{q\rightarrow +\infty}\left(\frac{h}{h^*_q}\right)^q = \lim_{q\rightarrow +\infty}\exp\left[\D q\ln\left(\!\frac{h}{h^*_q}\right)\right] = 0^+, \\[0.2cm]
\mbox{For } \D h=h^*_q & : & \psi_q(h^*_q) = \frac{1}{2}\,1\hspace{-0.15cm}1_{[0,h^*_q]}(h^*_q)\varphi(h^*_q)=\frac{1}{2}\,\varphi(h^*_q) \underset{q \to +\infty}{\longrightarrow} 0,
\end{eqnarray*}
as $h^*_q$ goes to $+\infty$ when $q$ goes to $+\infty$, $\varphi$ being a function with compact support.
\item The sequence of functions $\D\left(\psi_{q}\right)_{q\in\N^*}$ is uniformly dominated for all $q\in\N^*$ by an integrable function:
$$
\forall q\in N^*: |\psi_q(h)|\leq |\varphi|,
$$
and $|\varphi|\in L^1(\R)$ as $\varphi\in {\cal D(\R)}$.
\end{itemize}
The dominated convergence theorem enables us to conclude that
$$
\D  \lim_{q\rightarrow +\infty}\int_{-\infty}^{+\infty}\psi_q(h)dh \,= \int_{-\infty}^{+\infty}\lim_{q\rightarrow +\infty}\psi_q(h)dh \,= \int_{-\infty}^{+\infty}(H\varphi)(h)dh.
$$
With the same arguments, one gets for the second integral of (\ref{limit_distrib})
$$
\D  \lim_{q\rightarrow +\infty} \int_{-\infty}^{+\infty}\frac{1}{2}\!\left(\!\frac{h^{*}_q}{\!\!h}\!\right)^{\!q}1\hspace{-0.15cm}1_{[h^*_q,+\infty[}(h)\varphi(h)dh = 0,
$$
so that
$$
\D  \lim_{q\rightarrow +\infty} <T_{\mathcal{P}_{q}}, \varphi > \,\,= \int_{-\infty}^{+\infty}(H\varphi)(h)dh \,\, = \,\, <T_H, \varphi >, \forall \varphi \in {\cal D(\R)}.
$$
This ends the proof.
\end{prooff}
\noindent  In this setting, it is worth giving an interpretation of the results proved in this section. Basically, our results mean that when the distance between the values of $k_1$ and $k_2, (k_1 < k_2)$ increases, the finite elements  $P_{k_2}$ will be \emph{surely more accurate} than finite elements $P_{k_1}$, for values of the mesh size $h$ which fulfill the interval $[0,+\infty[$, and not only when $h$ goes to zero, as usually considered for accuracy comparison.\sa
Apart from the asymptotic case where $k_2-k_1$  goes to infinity, the probabilistic law (\ref{Nonlinear_Prob}) gives new insights into the relative accuracy between $P_{k_1}$ and $P_{k_2}$ finite elements. Indeed, for $h>h^*_{m,p}$, we obtained that $Prob\left\{ X_{m,p}^{(k_2)}(h) \leq X_{m,p}^{(k_1)}(h)\right\}\leq 0.5$. This shows that there are cases where $P_{k_2}$ finite elements \emph{probably} must be overqualified. As a consequence, a significant reduction of implementation time and execution cost could be obtained without loss of accuracy. Such a phenomenon has already been observed by using data-mining techniques coupled with other probabilistic models (see \cite{AsCh11, AsCh13, AsCh14}, \cite{AsCh16} and \cite{AsCh17}).
\section{Conclusion}\label{Conclusion}
\noindent In this paper, we derived an explicit $k-$dependence in $W^{m,p}$  \emph{a priori} error estimates, that we then applied to probabilistic relative accuracy of Lagrange finite elements. After having recalled some fundamental results of Banach spaces, especially the extension of C\'ea's classical Lemma to non Hilbert spaces, we derived general upper bounds on the basis functions and their partial derivatives for the polynomial space $P_k(K)$.\sa
Hence, we extended previous work \cite{ChasAs18}, \cite{Arxiv2} to the case of Banach $W^{m,p}$ spaces. This enabled us to evaluate  the relative accuracy between two Lagrange finite elements $P_{k_1}$ and $P_{k_2}, (k_1 < k_2)$, when the norm to measure the error estimate is defined on $W^{m,p}(\Omega)$. We also analyzed the asymptotic behavior of the relative accuracy between finite elements $P_{k_1}$ and $P_{k_1+q}$, for a fixed $k_1$, when $q$ goes to $+\infty$. We proved that, under some {\em ad hoc} assumptions that are fulfilled in most cases, the probabilistic law (\ref{Nonlinear_Prob}) is convergent to the Heaviside distribution $T_H$ in the weak-* topology on $D'(\R)$. \sa
Lastly, note that these perspectives are not  necessarily restricted to finite element methods, but can be extended to other approximation methods: given a class of numerical schemes and their corresponding error estimates, one can order them, not only by considering their asymptotic rates of convergence, but also by evaluating the most probably accurate one. \sa
\textbf{\underline{Homages}:} The authors want to warmly dedicate this research to pay homage to the memory of Professor Andr\'e Avez and Professor G\'erard Tronel, who largely promoted the passion of research and teaching in mathematics.


\begin{thebibliography}{100}
%
\bibitem{AgDn59} S. Agmon,  A. Douglis,  L. Nirenberg, Estimates near the boundary for solutions of elliptic partial differential equations satisfying general boundary conditions, {\em Comm. Pure Appl. Maths.} {\bf XII}, pp. 623--727 (1959).
%
\bibitem{Arcangeli_Gout} R. Arcangeli, J.L. Gout, Sur l'\'evaluation de l'erreur d'interpolation de Lagrange dans un ouvert de $\R^n$, \emph{ESAIM: Mathematical Modelling and Numerical Analysis - Mod\'elisation Math\'ematique et Analyse Num\'erique}, {\bf 10}, pp. 5--27 (1976).
%
\bibitem{AsCh11} F. Assous, J. Chaskalovic, Data mining techniques for scientific computing: Application to asymptotic paraxial approximations to model ultra-relativistic particles, \emph{J. Comput. Phys.}, {\bf 230}, pp. 4811--4827 (2011).
%
\bibitem{AsCh13} F. Assous, J. Chaskalovic, Error estimate evaluation in numerical approximations of partial differential equations: A pilot study using data mining methods, {\em C. R. Mecanique} \textbf{341} (2013) 304--313.
%
\bibitem{AsCh14} F. Assous, J. Chaskalovic, Indeterminate constants in numerical approximations of PDE's: A pilot study using data mining techniques, \emph{J. Comput. Appl. Math}, {\bf 270} (2014) 462-470.
\bibitem{BeSc72} M.S. Berger, M. Schechter, Embedding theorems and quasi-linear elliptic boundary value problems for unbounded domains, \emph{Trans. Amer. Maths. Soc.}, {\bf 172} (1972) 261-278.
%
\bibitem{BrSc08} S. C. Brenner, L. R. Scott, {\em The Mathematical Theory of Finite Element Methods}, Texts in Applied Mathematics, Vol. 15, 3rd ed., Springer, 2008.
%
\bibitem{Brezis} H. Brezis, Analyse fonctionnelle - Th\'eorie et applications, Masson (1992).
%
\bibitem{ChaskaPDE} J. Chaskalovic, Mathematical and numerical methods for partial differential equations, Springer Verlag, (2013).
%
%
%
\bibitem{AsCh16} J. Chaskalovic, F. Assous, Data mining and probabilistic models for error estimate analysis of finite element method, {\em Maths. And Comp. in Simulation} \textbf{129} (2016) 50--68.
%
\bibitem{AsCh17} J. Chaskalovic, F. Assous, Probabilistic approach to characterize quantitative uncertainty in numerical approximations, {\em Maths. Model. and Anal.} \textbf{22} (2013) 106--120.
%
\bibitem{ChasAs18} J. Chaskalovic, F. Assous, A new probabilistic interpretation of Bramble-Hilbert lemma,  {\em Comput. Meth.  Appl. Math.}, https://doi.org/10.1515/cmam-2018-0270 (2019).
%
\bibitem{Arxiv2} J. Chaskalovic, F. Assous, A new mixed functional-probabilistic approach for finite element accuracy, December 2018. arXiv:1803.09552 [math.NA]
%
\bibitem{Ciarlet_Raviart} Ciarlet, Raviart, General Lagrange and Hermite interpolation in $\R^n$ with applications to finite element methods, {\em Arch. Rat. Mech. Anal.} \textbf{46} (197) 177--199.
%
\bibitem{Ern_Guermond} A. Ern, J. L. Guermond, Theory and practice of finite elements, Springer, (2004).
%
\bibitem{OviFur} O. Furdui, Limits, Series, and Fractional Part Integrals, Springer, (2013).
%
\bibitem{Gris92} P. Grisvard, Singularities in boundary value problems, Springer, Berlin (1992).
%
\bibitem{Haubxx} Ch. Haubner, Finite Element Error Estimates in Non-Energy Norms for the Two-Dimensional Scalar Signorini Problem, {\em preprint}, (2019)
%
\bibitem{KuMi12} M. Kumar, G. Mishra, A Review on Nonlinear Elliptic Partial Differential Equations and Approaches for Solution,
{\em Intern. J. Nonlinear Sc.},  {\bf 13-4}, pp.401--418 (2012).
%
\bibitem{Lion69} J. L. Lions, Quelques m\'ethodes de r\'esolution des problemes aux limites non lin\'eaires, Dunod Gauthier-Villars, Paris, (1969).
%
\bibitem{Euler} E. Maor, "e:" The Story of a Number, Princeton Science Library (2015).
%
\bibitem{RaSc82} R. Rannacher, R. Scott, Some optimal error estimates for piecewise linear finite element approximations, {\em Math. Comp.}, {\bf 38} 437--445 (1982).
%
\bibitem{RaTho82} P.A. Raviart et J.M. Thomas, Introduction \`a l'analyse num\'erique des \'equations aux d\'eriv\'ees partielles, Masson (1982).
%
\bibitem{Schwartz} L. Schwartz, M\'ethodes math\'ematiques pour les sciences physiques, Hermann, (1983).
%
%
\end{thebibliography}
\end{document}